# A Theorem on Matroid Homomorphism

Jon Henry Sanders, jon_sanders@partech.com, ph. 315-335-7858

We call $f : M \to N$ a *matroid homomorphism*, or more briefly a *homomorphism* if $f : E(M) \to E(N)$ is an onto map of the ground set $E(M)$ of a matroid M onto the ground set $E(N)$ of a matroid N which preserves circuits of M, i.e., $f(C) \, \varepsilon \, \mathcal{C}(N)$ whenever $C \, \varepsilon \, \mathcal{C}(M)$, where $\mathcal{C}(X)$ denotes the circuit set of a matroid X. If $f^{-1}$ preserves circuits as well, we call f a *homeomorphism*. If M has no coloops and $f : M \to N$ is a homeomorphism it is not hard to show that M is isomorphic to a subdivision of N and that the elements of $f^{-1}(x)$ are in series for each $x \, \varepsilon \, E(N)$. We call a 1-1 homomorphism $f : M \to N$ a *circuit injection*. In this case N is isomorphic to a refinement of M. In this note we prove that when M is connected and N is binary and does not consist of a single circuit then any homomorphism $f : M \to N$ can be written as a composite,

$f = h \circ g$, of two homomorphisms h and g, where $g : M \to H$ is a homeomorphism and $h : H \to N$ is a circuit injection ( and thus M is (isomorphic to) a subdivision of H and N is (isomorphic to) a refinement of H. This implies that M must be binary. This theorem generalizes a result contained in a previous paper [ J. Sanders, *Circuit preserving edge maps II, J. Combin. Theory Ser. B* 42 (1987), 146-155] which derives the same decomposition under the (stronger) assumption that M and N are graphic.

## 1. Introduction

A matroid M *is of co-rank k* if $|E(M) - r(M)| = k$. We write " M is $CR^k$ " if M is connected and of co-rank k. The following facts from elementary set and matroid theory will be assumed.

<u>Fact 1.</u>  For any function  $f : S \to T$  from a set S onto a set T  we have for A, B $\subseteq$ S,

$f ( A \cup B) = f(A) \cup f(B)$,   $f(A) - f(B) \subseteq f(A - B)$, so  $f(A) \triangle f(B) \subseteq f(A \triangle B)$, where
A $\triangle$ B denotes the *symmetric difference* of two sets,  A $\triangle$ B  = (A-B) $\cup$ (B-A).

<u>Fact 2.</u>  Let M be connected, A $\subseteq$ E(M),  x $\in$ E(M), x $\notin$ A and M|A be $CR^k$ for some
positive integer k. Take B $\in$ $\mathcal{C}$(M) such that x $\in$ B, B meets A and B-A minimal.
Then M|(A $\cup$ B) is $CR^{k+1}$.

**Fact 3.** If M is $CR^2$ and A, B ∈ $\mathcal{C}$(M) then there exists a circuit C ∈ $\mathcal{C}$(M) such that
A Δ B ⊆ C.

**Fact 4.** If M|$E_1$, M|$E_2$ are $CR^2$ and M|($E_1 \cup E_2$) is $CR^3$ then $E_1 \cap E_2$ contains a circuit.

**Fact 5.** If M is $CR^2$ then M is a subdivision of the uniform matroid $U_{k,k+2}$ for some positive integer k.

**Fact 6.** If M is $CR^3$ and A, B, C ∈ $\mathcal{C}$(M) and B ≠ C and A∩B = ∅ then A∩C ≠ ∅.

**Fact 7.** If M is $CR^3$ and A, B ∈ $\mathcal{C}$(M), A ≠ B, A∩B ≠ ∅ and A ∪ B ⊂ E(M) then M| (A∪B) is $CR^2$.

## 2. Theorems

Let M be a connected matroid, f : M → B a homomorphism, B a binary matroid.

***Lemma 1.*** Let $x_1$, $x_2$ be distinct elements of E(M), f($x_1$) = f($x_2$), A ∈ $\mathcal{C}$(M), $x_1$ ∈ A, $x_2$ ∉ A. Let B ∈ $\mathcal{C}$(M) with $x_1$, $x_2$ ∈ B such that B−A is minimal. Then f(A) = f(B).

***Proof.*** Since f($x_2$) ∈ f(A) ∩ f(B), f($x_2$) ∉ f(A) Δ f(B). Also $x_2$ ∈ B −A so f($x_2$) ∈ f(B-A) ⊆ f (B-A) ∪ f(A-B) = f(AΔB). Thus by Fact 1, f(AΔB) properly contains f(A) Δ f(B). If f(A) ≠ f(B) then f(A) Δf(B) is a non-empty union of circuits. But M|(A∪B) is $CR^2$ (Fact 2 with k=1) so by Fact 3 there exists C ∈ $\mathcal{C}$(M|(A∪B)), C contains A Δ B so f(C) contains f(AΔB) which properly contains f(A) Δ f(B), a contradiction. Thus f(A) = f(B).

***Theorem 1.*** Let f : M → B be a homomorphism with M connected and B binary. Then either $f^{-1}$(x) is a series class for each x ∈ E(B) or B is a circuit.

***Proof.*** Assume $f^{-1}(x_0)$ is not a series class for some $x_0 \in E(B)$. Then there exist $x_1, x_2 \in f^{-1}(x_0)$, $A, B \in \mathcal{C}(M)$ which satisfy the hypothesis of Lemma 1. Thus $f(A \cup B) = f(A) = f(B) = D$ for some $D \in \mathcal{C}(B)$. If $E(M) = A \cup B$ we are through so assume there exists $c \in E(M)$, $c \notin A \cup B$. We will show $f(c) \in D$. (Since c was arbitrary this implies $f(M) = D$ and the theorem will be established). Since M is connected there exists $C \in \mathcal{C}(M)$, with $c \in C$ and $C \cap (A \cup B) \neq \emptyset$. Take C so that $C-(A \cup B)$ is minimal. Then $M|(A \cup B \cup C)$ is $CR^3$ by Fact 2 with k=2.

<u>Case 1.</u> There exists $x_1 \in D$ such that $f^{-1}(x_1) \cap (A \cup B)$ is not comparable to C, i.e., there exist $x, y \in f^{-1}(x_1) \cap (A \cup B)$, $y \in C$, $x \notin C$. Let $A' \in \mathcal{C}(M|(A \cup B))$, with $x, y \in A'$ and $A' - C$ minimal. Then $M|(C \cup A')$ is $CR^2$ and $x, y, C, A'$ satisfy the hypothesis of Lemma 1 so $f(C \cup A') = D'$ for some $D' \in \mathcal{C}(B)$. But Fact 4 then implies $D = D'$ so $f(c) \in D$.

<u>Case 2.</u> $f^{-1}(x) \cap (A \cup B)$ is comparable to C for each $x \in D$. Let $x_1 \in D$ be such that $f^{-1}(x_1) \cap C \cap (A \cup B)$ is empty. (If none such exists $f(C)$ contains D so $f(C) = D$ so we are done.) Let $P_1, ..., P_k$, (with k larger or equal to 3) be the maximal series classes of $M|(A \cup B)$. Then $\mathcal{C}(M|(A \cup B)) = \{(A \cup B) - P_j, j = 1, ..., k\}$ (Fact 5). Now there exist in $f^{-1}(x_1) \cap (A \cup B)$ distinct $x_i, x_j$, such that $x_i \in P_i \cap f^{-1}(x_1)$, $x_j \in P_j \cap f^{-1}(x_1)$. (Otherwise $f^{-1}(x_1) \cap (A \cup B) \subseteq P_h$ for some h and $C' = f((A \cup B) - P_h)$ is a circuit, $C' \subseteq D - x_1$ so $C' \subset D$, a contradiction.) Now of the sets $C \cap ((A \cup B) - P_i)$ and $C \cap ((A \cup B) - P_j)$ at least one must be non-empty (Fact 6), say the first. Also $x_i \notin ((A \cup B) - P_i) \cup C \subseteq A \cup B \cup C$ so $M|((A \cup B) - P_i) \cup C$ is $CR^2$ (Fact 7). Now $x_j \in ((A \cup B) - P_i) \triangle C$ and by Fact 3 there exists $C' \in \mathcal{C}(M)$, $C' \subseteq$

$((A \cup B) - P_i) \cup C$, $C' \supseteq ((A \cup B) - P_i) \Delta C$ so $x_i \notin C'$, $x_j x_j \in C'$ (and $A \cup B \cup C') = A \cup B \cup C$ ) and we have the conditions of Case 1 with $C'$ exchanged for C. Thus $f(c) \in D$.

**_Theorem 2._** Let $f : M \to N$ be a homomorphism such that $f^{-1}(x)$ is a series class for each $x \in E(N)$. Then $\{f(C)| C \in \mathcal{C}(M)\} := \mathcal{C}'$ is the set of circuits of a matroid H on E(N) and $f^{-1}(C') \in \mathcal{C}(M)$ when $C' \in \mathcal{C}' = \mathcal{C}(H)$, i.e., the function f, $f : E(M) \to E(N) = E(H)$ is a homeomorphism from M onto H which we denote by g, $g : M \to H$.

**_Proof._** Let $C' \in \mathcal{C}'$, $C' = \{x_1, ..., x_k\}$. We claim $C = f^{-1}(x_1) \cup ... \cup f^{-1}(x_k) = f^{-1}(C')$ is a circuit of M and also that $D \in \mathcal{C}(M)$, $f(D) = C'$ implies $D = C$. By definition there exists $X \in \mathcal{C}(M)$ such that $f(X) = C'$. Since, in general, $f^{-1}(f(X)) \supseteq X$ but for no circuit Y is $Y \subset C$ (for if such existed then $Y \cap f^{-1}(x_i)$ is empty for some i, $1 \leq i \leq k$ (since $f^{-1}(x_i)$ are in series) but then $f(Y) \subset C'$, a contradiction ) we have $C = X$. H is a matroid as follows. Let $y \in C'_1 \cap C'_2$, where $C'_1, C'_2 \in \mathcal{C}(H)$ and let $x \in f^{-1}(y)$. Then $x \in f^{-1}(C'_1) \cap f^{-1}(C'_2)$. Choose $C \subseteq f^{-1}(C'_1) \cup f^{-1}(C'_2)$ such that $x \notin C$. Then $f^{-1}(y) \cap C$ is empty and $f(C) \in \mathcal{C}(H)$, $f(C) \subseteq C'_1 \cup C'_2$, $y \notin f(C)$. Also if $C'_2 \subseteq C'_1$ then $f^{-1}(C'_2) \subseteq f^{-1}(C'_1)$ so $f^{-1}(C'_2) = f^{-1}(C'_1)$ so $C'_2 = C'_1$. Thus H is a matroid. It is clear that M is isomorphic to a subdivision of H and that N is a refinement of H. This establishes

**_Theorem 3._** Let $f : M \to B$ be a homomorphism, M connected, B binary and not a circuit. Then $f = h \circ g$, where g is a homeomorphism, h is a circuit injection (and M is isomorphic to a subdivision of g(M) ).

If M is not binary there exist $A, B \in \mathcal{C}(M)$ such that $M| (A \cup B)$ is $CR^2$ and thus equals a subdivision of the uniform matroid $U_{k,k+2}$, with $k \geq 2$. Thus there exists $C \in \mathcal{C}(M)$, such that A∆B ⊂ C. Thus M can have no binary refinement. Also if M is (isomorphic to) a subdivision of H then M is binary iff H is binary. This implies

**_Theorem 4._** Let f : M → B be a homomorphism, M connected, B binary and not a circuit. Then M is binary.

There is a result of Seymour [1] which characterizes when a binary matroid M can have three elements not all contained in some circuit. It has as a direct consequence the following

**_Corollary (to Seymour's theorem)._** Let M be a simple, vertically 4-connected binary matroid such that there exists a set of three elements S = {a,b,c}, S ⊆ E(M) with the property that no circuit of M contains S, i.e., for all C ∈ 𝒞(M), S ⊆ C is not true. Then M is graphic (and a, b, c correspond to three edges adjacent to a vertex).

This surprising result can be used with Theorem 4 above and the result of [2] to establish the following

**_Theorem 5._** Let f : M → G be a matroid homomorpism, M vertically 4-connected with no loops, G graphic and not a circuit. Then f is an isomorphism.